\documentstyle{amsppt}
\magnification=\magstep1
\NoBlackBoxes


\topmatter

\title A Banach subspace of $L_{1/2}$ which does not embed in $L_1$
(isometric version)
\endtitle

\author Alexander Koldobsky \endauthor
\address Division of Mathematics, Computer Science, and Statistics,
University
of Texas at San Antonio, San Antonio, TX 78249, U.S.A. \endaddress
\email koldobsk\@ringer.cs.utsa.edu \endemail

\abstract For every $n\geq 3,$ we construct an $n$-dimensional Banach
space which is isometric to a subspace of $L_{1/2}$ but is not
isometric
to a subspace of $L_1.$ The isomorphic version of this problem (posed
by S. Kwapien in 1969) is still open. Another example gives a Banach
subspace of  $L_{1/4}$ which does not embed isometrically in $L_{1/2}.$
Note that, from the isomorphic point of view, all the spaces $L_q$ with
$q<1$ have the same Banach subspaces.
\endabstract


\rightheadtext{A Banach subspace of $L_{1/2}$ which does not embed in
$L_1$}

\endtopmatter \document \baselineskip=14pt

\head 1. Introduction \endhead

A well known fact is that the space $L_1$ is isometric to a subspace
of $L_q$ for every $q<1.$ It is natural to ask whether the spaces $L_q$
with $q<1$ contain any Banach space structure not generated by $L_1.$
This question was first formulated in 1969 by Kwapien \cite{6} in
the following form: need every Banach subspace of $L_0$ be also a
subspace of $L_1?$ Later the question was mentioned by Maurey \cite{8,
Question 124}.

In 1970, Nikishin \cite{9} proved that every Banach subspace of $L_0$
is isomorphic
to a subspace of $L_q$ for every $q<1.$ Therefore, if we replace the
space $L_0$ in Kwapien's question by any of the spaces $L_q$ with $q<1$
we get an equivalent question.

Since all the spaces $L_q$ with $q<1$ embed in $L_0,$ Nikishin's
result also shows that these spaces are all the same
from the isomorphic Banach space point of view. Namely, every Banach
space which is isomorphic to a subspace of $L_q$ with $q<1$ is also
isomorphic to a subspace of $L_p$ for every other $p<1.$

In this paper we show that the answer to the isometric version of
Kwapien's question is negative. For every $n\in N,\ n\geq 3$
there exists an $n$-dimensional
Banach space which is isometric to a subspace of $L_{1/2}$ but is not
isometric to a subspace of $L_1.$ Using this example it is easy to see
that the spaces $L_q$ with $q<1$ may be different from the isometric
Banach space point of view. We give, however, a direct example
illustrating the difference by constructing a Banach subspace of
$L_{1/4}$ which does not embed isometrically in $L_{1/2}.$

The isomorphic version of Kwapien's question is still open. The most
recent related result seems to be a theorem of Kalton \cite{2} who
proved
that a Banach space $X$ embeds in $L_1$ if and only if $\ell_1(X)$
embeds
in $L_0.$

The isometric version of Kwapien's question can be reformulated
in the language of positive definite functions and stable random
vectors.
In fact, a Banach space $(X,\|\cdot\|)$ is isometric to a subspace
of $L_p$ with $0<p\leq 2$ if and only if the function $\exp(-\|x\|^p)$
is
positive definite \cite{1}. The main example of this paper
gives a norm such that the function $\exp(-\|x\|^{1/2})$ is positive
definite, but the function $\exp(-\|x\|)$ is not positive definite.
This result is close to problems of Schoenberg's type (see \cite{4}).

\head 2. The idea of the construction \endhead

Let $f$ be an infinitely differentiable even function on the unit
sphere
$S_{n}$ in $\Bbb R^n.$ We spoil the Euclidean norm $\|x\|_2$ in $\Bbb
R^n$ by means of the function $f.$ Namely, for $\lambda > 0$ consider
the function

$$\Cal N_{\lambda}(x) = \|x\|_2 \big( 1 + \lambda f({x\over {\|x\|_2}})
\big),\ x \in \Bbb R^n. \tag{1}$$

One can choose $\lambda$ small enough so that
$\Cal N_{\lambda}$ is a norm in $\Bbb R^n.$ This follows from a simple
one-dimensional
consideration: if $a,b\in R,$ $g$ is a convex function on $[a,b]$ with
$g''>\delta >0$ on $[a,b]$ and $h\in C^{2}[a,b]$ then the functions
$g + \lambda h$ have positive second derivatives on $[a,b]$ for
sufficiently small $\lambda$'s, and, hence, are convex on $[a,b].$

Let $\lambda_f = \sup \{\lambda > 0: \ \Cal N_{t}$ is a norm
in  $\Bbb R^n$ for every $t \leq \lambda \}.$ For each
$\lambda \leq \lambda_f,$ we denote by
$X_{\lambda}$ the Banach space with the norm
$\|x\|_{\lambda} = \Cal N_{\lambda}(x).$

Theorem 2 from the paper \cite{5} shows that, for every $q>0$ which is
not an even integer, there exists a small enough number $\lambda$ such
that the space $X_t$ is isometric to a subspace of $L_q$ for
every $t \leq \lambda.$
This fact was used
in \cite{5} to prove that, for every compact subset of $(0,\infty)
\setminus\{2k, k\in N\}$ there exists a Banach space different from
Hilbert spaces which is isometric to a subspace of $L_{q}$ for every
$q\in Q.$

For $q\in (0,1],$ let
$\lambda_q = \sup \{\lambda > 0: \ X_{t}$ is isometric to a
subspace of $L_q$ for every $t \leq \lambda \}.$  If
$0<q<p\leq 2$ then the space $L_p$
is isometric to a subspace fo $L_q,$ therefore, $\lambda_p \leq
\lambda_q.$
In particular,  $\lambda_1 \leq \lambda_q$ for every $q<1.$ Clearly,
$\lambda_1 \leq \lambda_f.$

Now we can explain the idea of getting a Banach subspace of $L_q$ with
$q<1$ which is not isometric to a subspace of $L_1.$ Suppose we can
find a function $f$ so that $\lambda_1$ is strictly less than
$\lambda_f$, and
also $\lambda_1$ is strictly less than $\lambda_q.$ Then, for every
$\lambda\in (\lambda_1,\ \min(\lambda_f,\lambda_q)],$ the space
$X_\lambda$
is a Banach space with the desired property.

Similarly, for $q<p<1,$ if we manage to find a function $f$ so that
$\lambda_p< \lambda_q$  and  $\lambda_p< \lambda_f$ we get an example
of a Banach space which embeds isometrically in $L_q$ but does
not embed in $L_p.$

The construction in \cite{5} is based on the use of spherical
harmonics and, in general, does not give a chance to calculate the
numbers $\lambda_q$ exactly. We are, however, able to choose a function
$f$ for which it is possible to calculate the numbers $\lambda_q$ for
certain values of $q.$ Our calculations do not depend on the results
from \cite{5} mentioned above, so the paper \cite{5} only shows a
direction for constructing
examples.
\bigbreak
We shall use one simple characterization of finite dimensional
subspaces
of $L_q.$

\proclaim{Proposition 1} Let $q$ be a positive number which is not an
even integer, $(X, \|\cdot\|)$ be an $n$-dimensional
Banach space, and suppose there exists a continuous function $b$ on
the sphere $S_n$ in $\Bbb R^n$ such that, for every $x\in \Bbb R^{n}$,
$$ \|x\|^{q} = \int_{S_{n}} |( x,\xi )|^{q}\ b(\xi)\ d\xi \tag{2}$$
where $(x,\xi)$ stands for the scalar product in $\Bbb R^n.$

Then $X$ is isometric to a subspace of $L_q$ if and only if
$b$ is a non-negative (not identically zero) function.
\endproclaim

\demo{Proof} If $b$ is a non-negative function we can assume without
loss of generality that $\int_{S_n} b(\xi)\ d\xi\ =\ 1.$ Choose any
measurable (with respect to Lebesgue measure) functions $f_1,\dots,f_n$
on $[0,1]$ having the joint distribution $b(\xi)d\xi.$ Then, by (2),
the operator $x\mapsto \sum x_i f_i,\ x\in \Bbb R^n$ is an isometry
from $X$ to $L_q([0,1]).$

Conversely, if $X$ is a subspace of $L_q([0,1])$ choose any functions
$f_{1},...,f_{n}\in L_q$ which form a basis in $X,$ and let $\mu$ be
the joint distribution of the functions $f_{1},...,f_{n}$ with respect
to Lebesgue measure. Then, for every $x\in R^{n}$,
$$\|x\|^{q}=\|\sum_{k=1}^{n} x_{k}f_{k}\|^{q}=\int_{0}^{1}
|\sum_{k=1}^{n}x_{k}f_{k}(t)|^{q} dt=\int_{R^{n}} |( x,\xi)|^{q}\ d\mu
(\xi)=\int_{S_{n}} |(x,\xi)|^{q}\ d\nu (\xi)$$
where $\nu$ is the projection of $\mu$ to the sphere.(For every
Borel subset $A$ of $S_{n}$,
$\nu(A)=\int_{\{tA,t\in R\}} \|x\|_{2}^{q} d\mu(x) ).$
It follows from (2) that
$$\int_{S_{n}} |( x,\xi )|^{q}\ b(\xi)\ d\xi =
\int_{S_{n}} |( x,\xi )|^{q}\  d\nu(\xi)$$
for every $x\in \Bbb R^n.$ Since $q$ is not an even integer,
we can apply the uniqueness theorem for measures on the sphere
from \cite{3} to show that $ d\nu(\xi) = b(\xi)\ d\xi$ which
means that $ b(\xi)\ d\xi$ is a measure, and the function $b$
is non-negative.
\qed \enddemo

The representation (2) exists for every smooth enough function
on the sphere (see, for example, Theorem 1 from \cite{5}). We are going
to choose special smooth norms
for which it is possible to calculate the function $b$ exactly
and then check if $b$ is non-negative. In this way we calculate
the numbers $\lambda_q$ for these norms.

We need the representation (2)
for some simple functions on the sphere.

\proclaim{Lemma 1} For every $x=(x_1,\dots,x_n)$ from the unit
sphere $S^n$ in $\Bbb R^n$ and every
$q>0$ we have
$$x_n^2 = {{\Gamma((n+q)/2)}\over{2\pi^{(n-1)/2} \Gamma((q+1)/2)}}
\int_{S_n} |(x,\xi)|^{q}\big( {{n+q}\over{q}}\xi_n^2 -
{1\over q}\big)\ d\xi, \tag{3}$$
Therefore,
$$x_n^2 = {{\Gamma((2n+1)/4)}\over{2\pi^{(n-1)/2}\Gamma(3/4)}}
\int_{S_n} |(x,\xi)|^{1/2}\big((2n+1)\xi_n^2- 2\big)\ d\xi, \tag{4}$$
and
$$x_n^2 = {{\Gamma((n+1)/2)}\over{2\pi^{(n-1)/2}}}
\int_{S_n} |(x,\xi)|\big((n+1)\xi_n^2- 1\big)\ d\xi. \tag{5}$$
Besides,
$$x_n^4 = {\Gamma((n+1)/2)\over{2\pi^{(n-1)/2}}}
\int_{S_n} |(x,\xi)|\big(-(n+3)(n+1)\xi_n^4 +
6(n+1)\xi_n^2- 3\big)\ d\xi. \tag{6}$$
\endproclaim

\demo{Proof} It is a well-known simple fact (see, for example,
\cite{7})
that for every $x\in \Bbb R^n$ and every $k>0,$
$$(x_1^2+\dots+x_n^2)^k = {{\Gamma((n+2k)/2)}\over{2\pi^{(n-1)/2}
\Gamma((2k+1)/2)}}\int_{S_n} |(x,\xi)|^{2k}\ d\xi. \tag{7}$$
Differentiate both sides of (7) by $x_n$ twice, and then use the fact
that $x\in S_n$ to get
$$1+(2k-2)x_n^2= (2k-1){{\Gamma((n+2k)/2)}\over{2\pi^{(n-1)/2}
\Gamma((2k+1)/2)}}\int_{S_n} |(x,\xi)|^{2k-2}\xi_n^2\ d\xi.$$
Use (7) with the exponent $2k-2$ instead of $2k$ to get
$$(2k-2)x_n^2= {1\over{2\pi^{(n-1)/2}}}\int_{S_n}
|(x,\xi)|^{2k-2}\Big((2k-1){{\Gamma((n+2k)/2)}
\over{\Gamma((2k+1)/2)}}\xi_n^2 - {{\Gamma((n+2k-2)/2)}
\over{\Gamma((2k-1)/2)}}\Big)\ d\xi.$$
Now use the fact that $\Gamma(x+1)=x\Gamma(x)$ and put $2k-2=q$
to get (3).

To prove (6), differentiate both sides of (7) four times by $x_n$
( remember that $x_1^2+\dots+x_n^2=1;$ do not factor the second and
the third derivatives !), and then put $k=5/2:$
$$-x_n^4 + 6x_n^2 + 3 = 4\Gamma((n+5)/2)\int_{S_n} |(x,\xi)| \xi_n^4
\ d\xi.$$
Now use (5), (7) with $k=1/2,$ and the fact that
$\Gamma(x+1)=x\Gamma(x)$
to get (6).
\qed \enddemo
\bigbreak

\head 3. Examples \endhead

For every $\lambda > 0$ define a function $\Cal N_{\lambda}$ on
$\Bbb R^n$ by

$$\Cal N_{\lambda}(x) = (x_1^2+\dots+x_n^2)^{1/2} \Big( 1 + \lambda
{{x_1^2+\dots+x_{n-1}^2-2x_n^2}\over {x_1^2+\dots+x_n^2}} \Big)^2,
\ x \in \Bbb R^n.$$
\medbreak
\proclaim{Lemma 2} $\Cal N_{\lambda}$ is a convex function if and only
if
$\lambda \leq {1\over{11}}.$ \endproclaim

\demo{Proof} The function $\Cal N_{\lambda}$ is convex
if and only if the following function of two variables is convex:
$$g(x,y) = (x^2+y^2)^{1/2}\Big(1+\lambda
{{x^2-2y^2}\over {x^2+y^2}} \Big)^2.$$
Calculating the second derivatives of the function $g$ we get
$$a^2{{\partial^2 g}\over{\partial x^2}} +
2ab {{\partial^2 g}\over{\partial x\partial y}}+
b^2 {{\partial^2 g}\over{\partial y^2}} = $$
$$(x^2+y^2)^{-7/2}(ay-bx)^2 \big(x^4(1-10\lambda-11\lambda^2)+
x^2 y^2(2-2\lambda+104\lambda^2) + y^4(1+8\lambda-20\lambda^2)\big).$$

The function $g$ is convex if and only if the latter expression is
non-negative for every choice of $a,b,x,y.$ Clearly, it happens
if and only if ${{-1}\over {10}}\leq \lambda \leq {{1}\over {11}}.$
\qed \enddemo
\bigbreak

For $\lambda \leq {{1}\over {11}}$ denote by $X_{\lambda}$ the Banach
space with the norm $\Cal N_{\lambda}.$

\proclaim{Theorem 1} Let $n\geq 3.$ If
$$\alpha_n = {{(18n^2-18n)^{1/2}-3n+1}\over{9n^2-12n-1}}\leq \lambda
\leq {{1}\over{6n-4}}$$
then the Banach space $X_{\lambda}$ is isometric
to a subspace
of $L_{1/2},$ and, at the same time, $X_{\lambda}$ is not isometric to
a subspace of $L_1.$
\endproclaim

\demo{Proof} Let us first prove that the space $X_\lambda$ is isometric
to a subspace of $L_{1/2}$ if and only if $ \lambda \leq
{{1}\over{6n-4}}.$
For every $x\in S_n,$ use (4) and (7) with $k=1/4$ to get the
representation (2) (with $q=1/2)$ for the norm $\Cal N_{\lambda} :$
$$\Cal N_\lambda^{1/2}(x) = 1 + \lambda (1-3x_n^2)= $$
$${{\Gamma((2n+1)/4)}\over{2\pi^{(n-1)/2}\Gamma(3/4)}}
\int_{S_n} |(x,\xi)|^{1/2}\big(1+ 7\lambda-
(6n+3)\lambda \xi_n^2 \big)\ d\xi.$$
Clearly, the function $b(\xi)= 1+ 7\lambda- (6n+3)\lambda \xi_n^2$
is non-negative on $S_n$ if and only if $\lambda \leq
{{1}\over{6n-4}},$
and the fact we need follows from Proposition 1.

Let us show that $X_\lambda$ is isometric to a subspace of $L_{1}$ if
and only if $ \lambda \leq \alpha_n.$ If $x\in S_n$ then
$\Cal N(x)= 1+2\lambda(1-3x_n^2)+ \lambda^2(1-6x_n^2+9x_n^4).$ We use
(5), (6) and (7) with $k=1/2$ to get the representation (2) for
$\Cal N_{\lambda}(x)$ with $q=1:$
$$\Cal N(x) = \Gamma((n+1)/2)\int_{S_n} |(x,\xi)|\ b(\xi_n^2) \ d\xi$$
where
$$b(\xi_n^2) = (1+8\lambda-20\lambda^2) +
\xi_n^2(n+1)(48\lambda^2-6\lambda)-
9\xi_n^4 (n+3)(n+1)\lambda^2.$$
The function $b$ is a
quadratic function of $\xi_n^2$ with negative first coefficient.
By Proposition 1, the space $X_\lambda$ embeds in $L_1$ if and only if
$b$ is non-negative for every $\xi_n^2\in [0,1].$ Clearly, it happens
if
and only if both numbers $b(0)=1+8\lambda-20\lambda^2$ and
$b(1)= 1+\lambda(2-6n)+\lambda^2(1+12n-9n^2)$ are non-negative.
Since $1+12n-9n^2$ is a negative number for every $n\geq 2$ and we
consider
only positive numbers $\lambda,$ it is clear that the condition is that
$ \lambda \leq \alpha_n.$
To prove the theorem, it suffices to note that, for every $n\geq 3,$
$1\over {6n-4}$ is less than $1\over {11},$  and $\alpha_n <
{1\over{6n-2}} < {1\over{6n-4}}.$ \qed \enddemo

\remark{Remarks} 1. Since every two-dimensional Banach space is
isometric to a subspace of $L_1$ it is impossible to construct
a two-dimensional space which the property of Theorem 1.
In our example, $\alpha_2 = {1\over{11}}$ which coincides with the
bound from Lemma 2.

2. The author is unable to apply the scheme from Section 2 to
every $q<1,$ although it is very likely the same idea works.
For instance, getting an example for $q=3/4$ is a matter of
calculating the eighth derivative of the function
$(x_1^2+\dots+x_n^2)^k.$
\endremark
\bigbreak

Theorem 1 shows also that the spaces $L_q$ with $q<1$ may have
different
Banach subspaces. The Banach subspaces of $L_{1/2}$ constructed in
Theorem 1 can not be isometric to subspaces of $L_q$ for all $q<1.$
In fact, if a Banach space $(X,\|\cdot\|)$ is isometric to a subspace
of $L_q$ for every $q<1$ then, by a theorem from \cite{1}, the function
$\exp(-\|x\|^q)$
is positive definite for every $q<1.$ The function $\exp(-\|x\|)$ is
then
positive definite as a pointwise limit of positive definite functions,
and the space $X$ (by the same result from \cite{1}) is isometric to a
subspace of $L_1.$
\bigbreak
On the other hand, we can show the difference directly.
For every $\lambda > 0$ define a function $\Cal N_{\lambda}$ on
$\Bbb R^3$ by

$$\Cal N_{\lambda}(x) = (x_1^2+x_2^2+x_3^2)^{1/2} \Big( 1 + \lambda
{{x_1^2+x_{2}^2-2x_3^2}\over {x_1^2+x_2^2+x_3^2}} \Big)^4,
\ x \in \Bbb R^3.$$
\medbreak

The proof of the following theorem is similar to the proofs of Lemma 2
and Theorem 1.
\proclaim{Theorem 2} The function $\Cal N_{\lambda}$ is a norm
if and only if $\lambda \leq {1\over{23}}.$
The corresponding Banach space $X_{\lambda}$ is isometric to
a subspace of $L_{1/4}$ if and only if $\lambda \leq {1\over{26}},$ and
$X_{\lambda}$ is not isometric to a subspace of $L_{1/2}$ if
$\lambda > {1\over{28}}.$ Thus, for ${1\over{28}} < \lambda
\leq {1\over{26}},$ the space $X_\lambda$ is a Banach subspace
of $L_{1/4}$ which does not embed isometrically in $L_{1/2}.$
\endproclaim

\enddocument